\documentclass{article}
\usepackage{amsmath}
\usepackage{amssymb}
\usepackage{amsfonts}
\usepackage{algorithm}
\usepackage[noend]{algpseudocode}

\newtheorem{theorem}{Theorem}
\newtheorem{acknowledgement}[theorem]{Acknowledgement}

\newtheorem{corollary}[theorem]{Corollary}

\newtheorem{lemma}[theorem]{Lemma}

\newtheorem{proposition}[theorem]{Proposition}
\newtheorem{remark}[theorem]{Remark}

\newenvironment{proof}[1][Proof]{\noindent\textbf{#1.} }{\ \rule{0.5em}{0.5em}}
\begin{document}

\title{Methods to construct the Sparse-paving Matroids over a Finite Set}
\author{B. Mederos$^{\ast \ast }$, M. Takane$^{\ast }$, G. Tapia-S\'{a}nchez$%
^{\ast \ast }$ and B. Zavala$^{\ast \ast \ast }$}
\maketitle

\begin{abstract}
In this work we present an algorithm to construct sparse-paving matroids over  finite set $S$. From this algorithm we derive some useful bounds on the cardinality of the set of circuits of any Sparse-Paving matroids  which allow us to prove in a simple way an asymptotic relation between the class of Sparse-paving matroids and the whole class of matroids. Additionally we introduce a matrix based method which render an explicit partition of the 
$r$-subsets of $S$, $\binom{S}{r}=\sqcup_{i=1}^{\gamma }\mathcal{U}_{i}$
such that each $\mathcal{U}_{i}$ defines a sparse-paving matroid of rank $r$.
 \end{abstract}
 
%\begin{keyword} 
Matroid, Paving matroid, Sparse-paving matroid,
Combinatorial Geometries, Lattice of a matroid.
%\end{keyword}

\section{Introduction}
We recall that a \textbf{matroid} $M=(S,\mathcal{I})$ consists of a finite
set $S$ and a collection $\mathcal{I}$ of subsets of $S$ (called the \textbf{
independent sets} of $M$) satisfying the following \textbf{independence
axioms}:
\begin{description}
\item{$\mathcal{I}1$} The empty set $\emptyset \in \mathcal{I}$.

\item{$\mathcal{I}2$} If $X\in \mathcal{I}$ and $Y\subseteq X$
then $Y\in \mathcal{I}$.

\item{$\mathcal{I}3$} Let $U,V\in \mathcal{I}$  with $
\left\vert U\right\vert =\left\vert V\right\vert +1$ then $\exists\, x \in
U\backslash V$ such that $V\cup \{x\}\in \mathcal{I}$.
\end{description}

A subset of $S$ which does not belong to $\mathcal{I}$ is called a \textit{
dependent set} of $M$. A \textbf{basis} [respectively, a \textbf{circuit}] of $M$ is a maximal
independent [resp. minimal dependent] set of $M$. A \textbf{basis [}
respectively, a \textbf{circuit] of }$M$ is a maximal [resp. minimal
dependent] set of $M$.\textbf{\ }The \textbf{rank} of a subset $X\subseteq S$
is $\text{rk}X:=\max \{\left\vert A\right\vert ;A\subseteq X$ and $A\in \mathcal{I}
\}$ and the rank of the matroid $M$ is $\text{rk}M:=\text{rk}S$. A \textbf{closed}
subset (or \textbf{flat}) of $M$ is a subset $X\subseteq S$ such that for
all $x\in S\backslash X$, $\text{rk}(X\cup \{x\})=\text{rk}X+1$. Then it can be defined
the \textbf{closure operator} $cl:\mathcal{P}S\rightarrow \mathcal{P}S$ on
the \textit{power set} of $S$, as follows: 
$$cl(X):=\min \{\,Y\subseteq
S\,:\,X\subseteq Y\; \text{and}\; Y\; \text{is closed in}\; M\,\}.$$
The \textbf{lattice of} a matroid 
$M$, denoted by $\mathcal{L}_{M}$ is the lattice defined by the closed sets
of $M$, ordered by inclusion where the meet is the intersection and the join
the closure of the union of sets. For general references of Theory of
Matroids, see \cite{W35}, \cite{NK2009}, \cite{We76} and \cite{O2011}. For
references of theory of lattices and theory of lattices of matroids, see 
\cite{Bi67}, \cite{H59}. \\

A matroid is \textbf{paving} if it has no circuits of cardinality less than
$\text{rk}M$. The matroid $M$ is named \textbf{sparse-paving} if $M$ and its dual $M^{\ast}$ are paving matroids.

In \cite{BCH73}, Blackburn, Crapo and Higgs asked: "In the enumeration of
(non-isomorphic) matroids on a set of $9$ or less elements, (sparse-)paving
matroids predominate. Does this hold in general?". There are several
results which suggest that the answer should be positive, see for example 
\cite{BPvP2012b}, \cite{BPvdP2012a}, \cite{Bo2010}. We give some methods to
construct the set of all sparse-paving matroids over a finite set $
S$ of any rank $r$, which allow us to give relations between the
cardinalities of $\left\vert Sparse_{n,r}\right\vert $ and $\left\vert
Matroid_{n,r}\right\vert$, where $\left( Sparse_{n,r}\right) Matroid_{n,r}$
is the set of all (sparse-paving) matroids of rank $r$ over a set $S$ of
cardinality $n$.

The manuscript is organized as follows: In Section 2, we give more
definitions and known results: For $n\geq 3$ and $rkM\geq 2$, there is an
equivalent definition of being a sparse-paving matroid. Namely, a paving
matroid $M=(S,\mathcal{I})$ with $\left\vert S\right\vert \geq 3$ and $\left\vert
\text{rk}M\right\vert \geq 2$ is a sparse-paving matroid if and only if its set of 
$\text{rk}M-$circuits $C_{rkM}$ satisfies the following property:
$$
\forall \ X,Y\in C_{rkM}
\Longrightarrow \left\vert X\cap Y\right\vert \leq  rkM-2.\ \ \ \ \ \ \ \
\ \ (**)
$$
Moreover, we give in lemma (1.3) the counterpart of the above result: Let 
$S$ be a set of cardinality $n$ and $2\leq r\leq n-1$. Then any set 
$C\subset \binom{S}{r}$ of $r-$subsets of $S$ satisfying property $(\ast
\ast)$ defines a matroid which is sparse-paving of rank $r$ with $C$ as its
set of $r-$circuits. The Section \ref{3s}  shows some bounds to the number of elements of the set of circuits of 
sparse-paving matroid as well as an algorithm to build sets with property (**).
Finally we end up this section giving another proof of the Piff's result \cite{P73}

$$\lim_{n\rightarrow \infty }\frac{\log _{2}\log _{2}\left\vert
Matroid_{n,\left [  \frac{n}{2} \right ]}\right\vert }{\log _{2}\log _{2}\left\vert
Sparse_{n,\left [  \frac{n}{2} \right ]}\right\vert }=1.$$

In Section \ref{4s}, we give an explicit construction of a partition of the 
$r-$subsets of $S$, $\binom{S}{r}=\sqcup_{i=1}^{\gamma }\mathcal{U}_{i}$
such that each $\mathcal{U}_{i}$ define a sparse-paving matroid of rank $r$
and $\gamma =2\max \left\{ \binom{r}{\left[ r/2\right] },\binom{n-r}{\left[
(n-r)/2\right] }\right\} $

\begin{acknowledgement}
The second author wants to thank Gilberto Calvillo for introduce her
into the theory of matroids, to Jes\'us de Loera and his group in the
University of California in Davis: this work was doing mainly during her
sabbatical-semester there; and to Criel Merino for useful discussions. This
research was partially supported by DGAPA-sabbatical-fellowship of UNAM. And
by Papiit-project $IN115414$ of UNAM.
\end{acknowledgement}

\section{A description of the Sparse-paving Matroids
through their set of circuits.}

\medskip The study of sparse-paving and paving matroids helps to understand
the behavior of the matroids in general and important examples of matroids
are indeed sparse-paving matroids, as the combinatorial finite geometries.
In 1959, Hartmanis \cite{H59} introduced the definition of paving matroid
through the concept of $d$-partition in number theory. Then, following a Rota\'
s suggestion, Welsh \cite[(1976)]{We76} called them, paving matroids. Later,
Oxley \cite{O2011} generalized the definition of paving matroid to include
all possible ranks. And Jerrum \cite{J2006} introduced the notion of
sparse-paving matroids.

\subsection{More definitions, notations and known results.}
Given a set $S$ of $n$ elements the following properties hold for the matroids with
base set $S$

\begin{description}
\item{a.} Any matroid $M=(S,\mathcal{I})$ is completely determined by itsf
set of basis, $\mathcal{B}$ which in turn determine the set $\mathcal{I}
=\{\,X\subseteq S\,:\,\exists B\in \mathcal{B} \;\text{with}\; X\subseteq B\,\}$.

\item{b.} Let $M=(S,\mathcal{I})$ be a matroid of rank rk$M$. Any circuit $X$ of $M$ has cardinality $\left\vert X\right\vert \leq \text{rk}M+1$.

\item{c.} Let $M=(S,\mathcal{I})$ be a matroid of rank $\text{rk}M$. Denote by $
M^{\ast }=(S,\mathcal{I}^{\ast})$ the \textbf{dual matroid} of $M$\textit{\ 
}whose set of basis is $\mathcal{B}^{\ast }:=S\backslash \mathcal{B} =\{S\setminus B\,:\,B \in \mathcal{B}\}$.
\end{description}

A matroid $M$ is called a \textbf{sparse-paving matroid} if $M$ and its dual 
$M^{\ast }$ are paving matroids.

\subsection{Examples of sparse-paving matroids}
In this section we present some well-known properties  
of the sparse-paving matroids as well as some examples and a new characterization
of this class of matroids. 

\begin{description}
\item{1.} A matroid is called \textbf{uniform} of rank $r$
over a set of $n$ elements, denoted by $U_{r,n}$, if all subsets of $S$ with cardinality $r$  are
basis. The dual matroid $U_{r,n}^{\ast }$ of a uniform matroid is again
uniform with rank $n-r$. Then any uniform matroid is sparse-paving. 

\item{2.} Any matroid of rank $1$ is paving since the empty
set is always independent.

\item{3.} By 1 and 2, for $n=1,2$ any matroid on $S$
is sparse-paving.
\end{description}

Therefore, along this paper we can assume that
$n\geq 3$ and $r\geq 2.$ Let $S$ be a set of cardinality $
\left\vert S\right\vert =n$. Denote by $\binom{S}{t}:=\{X\subseteq S  ;
\left\vert X\right\vert =t\}$ for $0\leq t\leq n$ the subsets of $S$ of
cardinality $t$ (called $t$-\textbf{subsets}). Let $M=(S,\mathcal{I})$ be a
\textbf{\ paving matroid} \textbf{of rank} rk$M$. Denote by $\mathcal{B}$
[resp. $\mathcal{C}_{rkM}$ ] the set of the basis [resp. rk$M$-circuits] of $M$. There exist an equivalent definition of being a sparse-paving matroid
given in \textbf{\cite{K74}\cite{Bo2010}}

\begin{lemma} 
\label{lemma1}
A paving matroid  $M=(S,\mathcal{I})$ with $\left\vert S\right\vert \geq 3$ and $\text{rk}M\geq 2$ is a sparse-paving matroid if and only if its set of $\text{rk}M$-circuits, $C_{\text{rk}M}$ satisfies the following property:
$$
\forall \ X,Y\in C_{rkM}\textit{\ with }X\neq Y \textit{we
have } \left\vert X\cap Y\right\vert \leq rkM-2\ \ \ \ \ \ \ \ \ \ (\ast
\ast )
$$
\end{lemma}

The next result is the counterpart of lemma \ref{lemma1}. That is, let $S$ be a set of cardinality $n\geq 3$ and $2\leq r\leq
n-1$. Then \textit{any} set $\mathcal{C\subseteq }\binom{S}{r}$ of $r$%
-subsets of $S$ satisfying property $(\ast \ast )$ defines a sparse-paving
matroid of rank $r$ with $\mathcal{C}$ as its set of $r$-circuits. In other
words, in this case, the ordered pair $(S,\mathcal{I})$ with $\mathcal{I}
:=\{X\subseteq S;$ $\exists B\in \binom{S}{r}\backslash \mathcal{C}\}$ is in
fact a matroid (ie., $\mathcal{I}$ satisfies the independent axioms of a
matroid, see Introduction) and paving. Then by lemma \ref{lemma1}, $(S,\mathcal{I})$ is sparse-paving.

\begin{proposition}
\label{prop}
Let $S$ be a set
of cardinality $\left\vert S\right\vert =n\geq 3$ and $2\leq
r\leq n-1.$ Let $C\subset \binom{S}{r}$ be a set of $r$-subsets of $S$, satisfying the following property

$$\forall X,Y\in C \textit{ with } X\neq Y \textit{ then }\left\vert X\cap
Y\right\vert \leq r-2.\ \ \ \ \ \ \ \ \ (\ast \ast )$$
\noindent
Define $M:=(S,\mathcal{I})$ where $\mathcal{B}:=\binom{S}{r}
\backslash C$ and 
$$\mathcal{I}:=\{ X\subseteq S\,:\,
\exists B \in \mathcal{B}\; \text{with}\; X\subseteq B\}.$$ 
Then, \textbf{
(A)} $M$ is a matroid of $rkM=r$ and \textbf{(B)} $M$ is sparse-paving.
\end{proposition}

\begin{proof}
 Let $S$ be a set and take a subset $\mathcal{C}\subset \binom{S}{r}$ 
 satisfying the property $(\ast \ast )$. Take $M=(S,
\mathcal{I})$ with set of basis $\mathcal{B}=\tbinom{S}{r}\backslash 
\mathcal{C}$. Let start by proving {\bf A}. Let $M$ be a matroid of rank $r$. For this proof, we will use an equivalent definition of matroid, which says:

\textit{Let $M=(S,\mathcal{I})$ is a matroid if and only if $\mathcal{I}$ satisfies (
$\mathcal{I}$1),($\mathcal{I}$2) as in the introduction and ($\mathcal{I}$3)$
^{\prime }$: let $B_{1},B_{2}\in \mathcal{B}$ be two basis of $M$ and $x\in
B_{1}\backslash B_{2}.$ To prove $\exists y\in B_{2}\backslash B_{1}$ such
that $(B_{1}\backslash \{x\})\cup \{y\}\in \mathcal{B}$}. \\
\noindent
Two case are possible 
\begin{enumerate}
\item  $\left\vert S\right\vert =3$ and rk$M=2,$ the
possibilities for $\mathcal{C}$\ to have property $(\ast \ast )$ are $
\mathcal{C}=\emptyset $ or $\left\vert \mathcal{C}\right\vert =1.$ In both
cases, $M$ is a matroid and it is sparse-paving.

\item $\left\vert S\right\vert \geq 4.$
\begin{description}
\item [Proof of $(\mathcal{I}1)$].  It is enough to prove that $\mathcal{B}$ is not empty. Since $n\geq 4$, $2\leq r\leq n-1$ and $S=\{1,...,r,r+1,...,n\}.$ Take $%
A_{1}=\{1,...,r-1,r\},$ $A_{2}=\{1,...,r-1,r+1\}$ which are subsets of $S$
with cardinality $r$ and $\left\vert A_{1}\cap A_{2}\right\vert =r-1$. Then
by $(\ast \ast )$, there exists $ i\in \{1,2\}$ such that $A_{i}\in \mathcal{B}$.
Then $\mathcal{B}\neq \emptyset$.

\item [Proof of $(\mathcal{I}2)$]. Let $Y\subseteq X\subseteq S$ such that there exists $B\in \mathcal{B}$ with $X\subseteq B$. Then $Y\subseteq B$, that is 
$Y$ is independent by definition.

\item [Proof of $(\mathcal{I}^{\prime }3)$].  let $B_{1},B_{2}\in 
\mathcal{B}$ be two basis of $M$ and $x\in B_{1}\backslash B_{2}$. To prove $
\exists y\in B_{2}\backslash B_{1}$ such that $(B_{1}\backslash \{x\})\cup
\{y\}\in \mathcal{B}$ we deal with two cases:

\begin{itemize}
\item Assume $m:=\left\vert
B_{2}\backslash B_{1}\right\vert =1$. That is, $B_{2}\cap
B_{1}=B_{1}\backslash \{x\}$ and $B_{2}=(B_{1}\backslash \{x\})\cup \{y\}$
for some $y\in S.$ Then $(B_{1}\backslash \{x\})\cup \{y\}\in \mathcal{B}$.

\item Let define $m:=\left\vert
B_{2}\backslash B_{1}\right\vert \geq 2$ and let 
$$B_{2}=(B_{1}\cap
B_{2})\cup \{y_{1},y_{2},y_{3},...,y_{m}\}.$$
Define $A_{i}:=(B_{1}\backslash
\{x\})\cup \{y_{i}\}$ for $i=1,...,m.$ Since $\forall i\neq j,$ $\left\vert
A_{i}\cap A_{j}\right\vert =r-1$ and $m\geq 2$, by $(\ast \ast ),$ $\exists
A_{i_{0}}\in \mathcal{B}$. Therefore, $(B_{1}\backslash \{x\})\cup
\{y_{_{i_{0}}}\}=A_{i_{0}}\in \mathcal{B}$, and $M$ \textit{is a matroid}.
\end{itemize}
\end{description}
\end{enumerate}

Next we prove $\mathbf B$. First we will prove that $M$ is a \textit{paving
matroid}, which is equivalently to prove $\forall Z\subseteq S$ of $\left\vert
Z\right\vert =rkM-1,$ $Z\in \mathcal{I}$. This proof is similar to the
one of ($\mathcal{I}1$). Namely:

Let consider rk$M\leq n-1$. Since $n\geq 3$ and $\left\vert Z\right\vert =rkM-1$ we have $S=Z\cup \{x_{1},x_{2},...,x_{m}\}$ with $m\geq 2$. Let denote $
A_{i}:=Z\cup \{x_{i}\}$ for $i=1,2,...,m$. By $(\ast \ast )$ and $m\geq 2, 
\exists i_{0}\in \{1,...,m\}$ such that $Z\subset A_{i_{0}}\in \mathcal{B}$.
Then $Z\in \mathcal{I}$. To conclude the proof we use lemma 2 and obtain that $M$ is a sparse-paving matroid.
\end{proof}

\section{ A bound for the cardinality of the set of 
$\text{rk}M$ circuits of a sparse-paving matroid $M$}
\label{3s}

In this section we give an easy-finding bounds of the set of $r$-circuits of
the sparse-paving matroids. Recall from lemma \ref{prop}, the  property $(\ast \ast
)$ on a set $\mathcal{U}$  to be the collection of $r$-circuits of a
matroid: For all $X,Y\in \mathcal{U}$ we have $\left\vert X\cap Y\right\vert
\leq r-2$.

\subsection{Technical steps to construct sets with property $(\ast
\ast )$}

In this section before to state the main results we start by giving some
technical observations. \\

Let $S$ with $\left\vert S\right\vert =n$, $2\leq r\leq n-1$ and $%
X_{1}\in \binom{S}{r}$ be fixed. To find the set $\{X\in \binom{S}{r}%
;\left\vert X\cap X_{1}\right\vert \leq r-2\}$ is enough to find the set $%
\{A\in \binom{S}{r};\left\vert A\cap X_{1}\right\vert =r-1\}.$ One can
readily deduce the following remarks

\begin{remark} 
\label{1} $$\{A\in \binom{S}{r};\left\vert X_{1}\cap
A\right\vert \leq r-2\}=\binom{S}{r}\backslash \{A\in \binom{S}{r}%
;\left\vert X_{1}\cap A\right\vert =r-1\}.$$
\end{remark}
\begin{remark}
\label{2} $\left\vert \{A\in \binom{S}{r};\left\vert X_{1}\cap
A\right\vert =r-1\}\cup \{X_{1}\}\right\vert =r(n-r)+1$. Let us take all the $r+1$-subsets of $S$ containing $X_{1}$

$$\{Y_{1,1},Y_{1,2},...,Y_{1,n-r}\}=\{Y\in \binom{S}{r+1}\,:\,
X_{1}\subset Y\}.$$
\end{remark}
That is, $Y_{1,i}=X_{1}\cup \{v_{i}\}$ where $S\backslash
X_{1}=\{v_{1},...,v_{n-r}\}$. Let $X_{2}\in \{A\in \binom{S}{r};\left\vert
X_{1}\cap A\right\vert \leq r-2\}$ we can construct $
\{Y_{2,1},Y_{2,2},...,Y_{2,n-r}\}=\{Y\in \binom{S}{r+1};X_{2}\subset Y\}.$
The following observation holds

\begin{remark} 
\label{3}
Let us take all the $r+1$-subsets of $S$ containing $X_1$
$$\{Y_{1,1},Y_{1,2},...,Y_{1,r-n}\}.$$
The following observation holds.
If $\left\vert X_{1}\cap X_{2}\right\vert
\leq r-2$ then 
$$\{Y_{1,1},Y_{1,2},...,Y_{1,r-n}\}\cap
\{Y_{2,1},Y_{2,2},...,Y_{2,r-n}\}=\emptyset.$$
\end{remark}
\begin{proof} Note that $X_{1}\nsubseteqq Y_{1,k}$ and $X_{2}\subset
Y_{1,k}$ for $k=1,..,n-r.$ Similarly $X_{2}\nsubseteqq Y_{2,j}$ and $%
X_{1}\subset Y_{2,j}$ for $j=1,..,n-r.$ Therefore $Y_{1,k}\neq Y_{2,j}$ for $%
k,j\in \{1,...,n-r\}$.
\end{proof}\\

\noindent
The  remark \ref{1} allows us to build the
following algorithm 

\begin{algorithm}[h]
\caption{Algorithm}\label{euclid}
\begin{algorithmic}[h]
\State \textbf{Set} $\mathcal{U} \gets\emptyset $
\State \textbf{Set }$B \gets\emptyset $
\State \textbf{Set }$i \gets 1$
\State \textbf{Repeat}
\State \;\; \textbf{Select} $X_{i}\in \binom{S}{r}/B$
\State \;\; \textbf{Do} $\mathcal{U} \gets \mathcal{U}\cup \{X_{i}\}$
\State \;\; \textbf{Do} $B_{i} \gets \{A\in \binom{S}{r};\left\vert X_{i}\cap A\right\vert
=r-1\}\cup \{X_{i}\}$
\State \;\;\;\textbf{Do} $B \gets B\cup B_{\substack{ i }}$
\State \;\; $i \gets i+1$.
\State \textbf{Until} $\binom{S}{r}/B=\emptyset $
\State \textbf{Ouput} $\mathcal{U}$
\end{algorithmic}
\end{algorithm}
The set $\mathcal{U}$ obtained from the above algorithm satisfies some interesting properties
stated in the following two lemmas

\begin{lemma}Given any set\textbf{\ }$\mathcal{U}$ obtained
from the above algorithm. Then $\mathcal{U}$ is a maximal set fulfilling 
the property $(\ast \ast )$.
\end{lemma}

\begin{proof} Let$X_{i},X_{j}\in \mathcal{U}$ with $i<j,$ then
by the algorithm $X_{j}\in \binom{S}{r}/B_{i}$, therefore $\left\vert
X_{i}\cap X_{j}\right\vert \leq r-2$ and consequently $\mathcal{U}$ satisfy
the property $(\ast \ast ).$ In order to prove that $\mathcal{U}$ is a
maximal set we suppose that there exists $C\in \binom{S}{r}$ satisfying $%
C\notin U$ and $\left\vert X_{i}\cap C\right\vert \leq r-2,$ $\forall
X_{i}\in U$. Then $C\notin B$ which is a contradiction with the fact that $B=%
\binom{S}{r}$.
\end{proof}\\

Base on remarks \ref{2} and \ref{3} we can find upper and lower bounds on the
maximal sets $\mathcal{U}$ that satisfy the property  $(\ast \ast )$.

\begin{lemma}
\label{bounds}
Let $S$ be a set of cardinality 
$n$ and let $2\leq r\leq n-1$.

\begin{description}
\item{a).} Assume that $\mathcal{U}\subset \binom{S}{r}$ satisfies
property $(\ast \ast )$.  Then $\left\vert \mathcal{U}\right\vert
\leq \frac{1}{n-r}\binom{n}{r+1}$.

\item{b).} For all maximal set $\mathcal{U}\subset \binom{S}{r}$
satisfying property $(\ast \ast )$. Then $\frac{1}{r(n-r)+1}%
\binom{n}{r}\leq \left\vert \mathcal{U}\right\vert$.
\end{description}
\end{lemma}

\begin{proof}
In order to prove the item a) let us consider 
$\mathcal{U=\{}X_{1,}X_{2},...,X_{k}\}$. For each $X_{i}\in \mathcal{U}$ we
construct the set $A_{i}=\{Y_{i,1},Y_{i,2},...,Y_{i,n-r}\}=\{Y\in \binom{S}{
r+1};$ $X_{i}\subset Y\}$. By remark \ref{3} we have that  $A_{l}\cap
A_{j}=\emptyset,\, l\neq j$. On the other hand, since $A=\underset{i=1}{
\overset{k}{\sqcup }}A_{i}$ and $A\subset \binom{S}{r+1}$ we get 

\begin{center}
$k(n-r)=\overset{k}{\underset{i=1}{\sum }}\left\vert A_{i}\right\vert
=\left\vert A\right\vert \leq \left\vert \binom{S}{r+1}\right\vert $
\end{center}
Therefore $\left\vert \mathcal{U}\right\vert \leq \frac{1}{n-r}
\binom{n}{r+1}$ concluding the prove of item a).\\

To carry out the prove of item b) we have to deal with the set 
$$B_{i}=\{A\in 
\binom{S}{r};\left\vert X_{i}\cap A\right\vert =r-1\}\cup \{X_{i}\}.$$
 Since the set $\mathcal{U=\{}X_{1,}X_{2},...,X_{k}\}$ is maximal, then $B=\underset
{i=1}{\overset{k}{\sqcup }}B_{i}=\binom{S}{r+1}$. Combining this with the
observation 2 we deduce 

$$k(r(n-r)+1)=\overset{k}{\underset{i=1}{\sum }}\left\vert B_{i}\right\vert
\geq \left\vert B\right\vert =\left\vert \binom{S}{r}\right\vert $$
from this inequality we obtain that $\frac{1}{r(n-r)+1}\binom{n}{r}\leq
\left\vert \mathcal{U}\right\vert $.
\end{proof}\\

From previous lemma we can bound the number of elements of the set
of circuits of any sparse-paving matroid, this is stated in the next
Corollary. For a similar result, see also \cite[(4.8)]{MNRV2010}.

\begin{corollary}  Let $S$ be a set of
cardinality $n$ and let $2\leq r\leq n-1$.  Let $(S,C_{n,r})$ be a sparse-paving matroid of rank $r$ and $
C_{n,r}$ its set of $r$-circuits. Then $$\left\vert 
\mathcal{C}_{n,r}\right\vert \leq \frac{1}{n-r}\binom{n}{r+1}.$$
\end{corollary}
\begin{proof} This follows directly from  lemma \ref{bounds}, item a) and
the lemma 2.
\end{proof}\\

Another consequence of the  lemma \ref{bounds} is that
it is possible to find a lower bound on the number of sparse-paving matroids

\begin{corollary}
\label{Secondbound} Let $S$ be a set of
cardinality $n$ and let $2\leq r\leq n-1$. Then

$$2^{\left[ \frac{1}{r(n-r)+1}\binom{n}{r}\right] }\leq \left\vert
Sparse_{n,r}\right\vert .$$
\end{corollary}

\begin{proof} By lemma \ref{bounds}, any $\mathcal{U}\subset 
\binom{S}{r}$ with property $(\ast \ast )$ satisfies $\frac{1}{r(n-r)+1}
\binom{n}{r}\leq \left\vert \mathcal{U}\right\vert $. Let $\mathcal{P(U)}
=\{X\subseteq \mathcal{U}\}$ be the power set of $\mathcal{U}$. Then $
\forall \mathcal{C}\in \mathcal{P(U)}$, $\mathcal{C}$ satisfies property $
(\ast \ast )$, and therefore $\mathcal{C}$ defines a sparse-paving matroid. Moreover,
if $\mathcal{C\neq C}^{\prime }$ in $\mathcal{P(U)}$, their respective
sparse-paving matroids are different. Since $\frac{1}{r(n-r)+1}\binom{n}{r}
\leq \left\vert \mathcal{U}\right\vert$, we have 
$$2^{\left[ \frac{1}{
r(n-r)+1}\binom{n}{r}\right] }\leq \left\vert \mathcal{P(U)}\right\vert
=2^{\left\vert \mathcal{U}\right\vert }\leq \left\vert
Sparse_{n,r}\right\vert.$$
\end{proof}\\

\noindent
From this lower bound we get the same result of Piff \cite{P73}.
\begin{corollary}
Let $S$ be a set of
cardinality $n$. Then

$$\lim_{n\rightarrow \infty }\frac{\log _{2}\log _{2}\left\vert
Matroid_{n, \left[ \frac{n}{2}\right]}\right\vert }{\log _{2}\log _{2}\left\vert
Sparse_{n,\left[ \frac{n}{2} \right]}\right\vert }=1.$$
\end{corollary}
\begin{proof}
It is readily  to see that  $\left\vert Matroid_{n, \left[ \frac{n}{2}\right]}\right\vert \le 2^{2^n}$. On the other hand, by corollary \ref{Secondbound} we have for $r=\left[ \frac{n}{2}\right]$ the following inequality
\begin{equation}
\label{lower}
2^{ \frac{1}{\left[ \frac{n}{2}\right](n-\left[ \frac{n}{2}\right])+1}\binom{n}{\left[ \frac{n}{2}\right]} }\leq \left\vert
Sparse_{n,\left[ \frac{n}{2}\right]}\right\vert .
\end{equation}
In order to easy the notation lets define $p:=\left[ \frac{n}{2}\right](n-\left[ \frac{n}{2}\right])+1$ and $p \sim  n^2$. Using the well-known Stirling's approximation for the binomial coefficients 
$$
\sqrt{\left[ \frac{n}{2}\right]}\binom{n}{\left[ \frac{n}{2}\right]} \ge 2^{2 \left[ \frac{n}{2}\right] -1}
$$
and substituting in \eqref{lower} we get 
\begin{equation}
\label{lower}
2^{ \frac{1}{p\sqrt{\left[ \frac{n}{2}\right]}} 2^{2 \left[ \frac{n}{2}\right] -1} }\leq \left\vert
Sparse_{n,\left[ \frac{n}{2}\right]}\right\vert .
\end{equation}
Consequently we have 
\begin{equation}
\label{ine}
1 \le \frac{\log _{2}\log _{2}\left\vert
Matroid_{n, \left[ \frac{n}{2}\right]}\right\vert }{\log _{2}\log _{2}\left\vert
Sparse_{n,\left[ \frac{n}{2} \right]}\right\vert } \le \frac{n}{2 \left[ \frac{n}{2}\right] -1 - \log_2\log_2{p\sqrt{\left[ \frac{n}{2}\right]}}}
\end{equation}
Note that $\log_2\log_2{p\sqrt{\left[ \frac{n}{2} \right]}}$ has a smaller order of growth than $n$. Therefore, applying limits to both sides of \eqref{ine}  the conclusion follows.
\end{proof}

In the next section will find lowers and uppers bounds for the numbers of $r$-sets satisfying the property $(\ast
\ast)$. Let us start by counting the number of sets 
$$\left\{ \left(
X_{1},...,X_{\alpha }\right) : \text{satisfying}\;(\ast \ast )\right\} $$
with $\alpha:=\min \left\{ \frac{1}{r}\binom{n}{r-1},\frac{1}{n-r}%
\binom{n}{r+1}\right\}$.

\subsection{Bounds for the cardinality of sets satisfying $(**)$  }
Now we will present some bounds related with  sets satisfying the condition $(**)$ aiming to obtain a lower and upper bound of $Sparse_{n,r}$  

\begin{remark}
	\label{r12} 
	Given $X\in \binom{S}{r}$. The  cardinality  of the set $\left\{ Y\in \binom{S}{r}:\left\vert Y\cap X\right\vert \le r-2\right\} $
	is equal to
	$$\binom{n}{r}-r(n-r)-1.$$	
\end{remark}

\begin{proof}We readily  see that 
	\begin{eqnarray}
	\left\{ Y\in \binom{S}{r}:\left\vert Y\cap
	X\right\vert \leq r-2\right\} &=& \binom{S}{r} \setminus \left ( \{X\}\cup \left \{ Y\in \binom{S}{r}
	:\left\vert Y\cap X\right\vert =r-1\right\} \right ) \nonumber \\
	&=& \binom{S}{r} \setminus \left (\{X\}\cup \left\{ Z\cup \{a\}\in \binom{S}{r}:Z\in \binom{X}{r-1},\;a\in S\backslash X\right\} \right )\nonumber
	\end{eqnarray}
	Therefore
	
	\begin{eqnarray}
	\left\vert \left\{ Y\in \binom{S}{r}:\left\vert
	Y\cap X\right\vert \le r-2\right\} \right\vert &=& \left\vert \binom{S}{r}
	\right\vert - \left\vert \binom{X}{r-1}
	\right\vert \left\vert S\backslash X\right\vert+1 
	\noindent \\ 
	&=&\binom{n}{r} - r(n-r)+1.
	\end{eqnarray}
\end{proof}\\	

\noindent
As a consequence of Remark \ref{r12}, the following is obtained

\begin{remark}
	\label{r13}	
	The cardinality  of the set $\left\{ \left( X_{1},X_{2}\right) :\text{ satisfying }
	(\ast \ast )\right\}$ is equal to
	$$\binom{n}{r}\left[ \binom{n}{r}-r(n-r)-1%
	\right].$$
\end{remark}
\begin{proof} This  is a direct consequence of  Remark 
	\ref{r12}.
\end{proof}\\

\noindent
Now we will establish an useful intersection lemma

\begin{lemma}
	\label{intersect}
	Given $\{X_1, X_2\}$ satisfying $(**)$. The set 
	$$\left\{ Y\in \binom{S}{r}
	:\left\vert Y\cap X_{1}\right\vert =r-1, \,\left\vert Y\cap
	X_{2}\right\vert =r-1\right\} \neq \emptyset.$$
if and only if
	$$\left\vert X_{1}\cap X_{2}\right\vert =r-2.$$
\end{lemma}

\begin{proof}
	Let us consider 
	$$\left\{ Y\in \binom{S}{r}
	:\left\vert Y\cap X_{1}\right\vert =r-1, \,\left\vert Y\cap
	X_{2}\right\vert =r-1\right\} \neq \emptyset.$$
	Then, there exists 
	$$T\in \left\{ Y\in \binom{S}{r}:\left\vert Y\cap X_{1}\right\vert=r-1\right\} \cap \left\{ Y\in \binom{S}{r}:\left\vert Y\cap X_{2}\right\vert =r-1\right\}.$$
	From this it is obtained that 
	\begin{equation}
	\label{T}
	T=\left( T\cap X_{1}\right) \cup a=\left( T\cap X_{2}\right) \cup b,
	\end{equation}  
	where
	\begin{equation}
	\label{17}
	a\notin X_{1}\; \text{and}\; b\notin X_{2}.
	\end{equation}
	We claim that $a\neq b$. Let us suppose that $a=b$. Then
	$$ T\cap X_{1}=T\cap X_{2},$$
	which implies that $T\cap X_{1}=T\cap X_{2}\subseteq X_{1}\cap X_{2}$, which is a contradiction due to the fact
	$\left\vert T\cap X_{1}\right\vert =r-1$ and $\left\vert X_{1}\cap X_{2}\right\vert \leq r-2.$\\
	\noindent
	Since $a\neq b$ from \eqref{T} we deduce that $a\in T\cap X_{2}$ and $b\in T\cap X_{1}$. The equation \eqref{T} and \eqref{17} imply that 
	\begin{eqnarray}
	\label{TT}
	T &= & \left\langle \left( T\cap X_{1} \right) \cup \{a\} \right\rangle \cap \left\langle \left( T\cap X_{2} \right ) \cup  \{b \} \right \rangle \nonumber \\
	& = & \left\langle \left( T\cap
	X_{1}\right) \bigcap \left[ \left( T\cap X_{2}\right) \cup \{b\}\right]
	\right\rangle \bigcup \left\langle \{a\}\bigcap \left[ \left( T\cap X_{2}\right)
	\cup \{b\}\right] \right\rangle \nonumber \\
	&=&\left[ \left( T\cap X_{1}\right) \cap \left( T\cap X_{2}\right) \right]
	\bigcup \left[ \left( T\cap X_{1}\right) \cap \{b\}\right] \bigcup \left[ \{a\}\cap
	\left( T\cap X_{2}\right) \right] \bigcup \left[ \{a\}\cap \{b\}\right] \nonumber \\
	&=&\left[ T\cap X_{1}\cap X_{2}\right] \cup \{b\}\cup \{a\}.	
\end{eqnarray}	
	From \eqref{TT} we obtain the following inequality
	
	$$r=\left\vert T\right\vert =\left\vert \left[ T\cap X_{1}\cap X_{2}\right] \cup \{b\}\cup \{a\}\right\vert =\left\vert T\cap X_{1}\cap X_{2}\right\vert +2.$$
	Therefore	
	$$r-2=\left\vert T\cap X_{1}\cap X_{2}\right\vert \leq \left\vert X_{1}\cap X_{2}\right\vert \leq r-2,$$
	which leads to
	$$\left\vert X_{1}\cap X_{2}\right\vert =r-2.$$
	\noindent
	Now we prove the converse implication. If $\left\vert X_{1}\cap X_{2}\right\vert =r-2$, then
	
	$$X_{1}=\left( X_{1}\cap X_{2}\right) \cup \{x,y\}$$ 
	and
	$$X_{2}=\left(X_{1}\cap X_{2}\right) \cup \{u,v\}$$
	with $x,y\notin X_{2}$ and $u,v\notin
	X_{1}.$
	Let us denote 
	$$A=\left\{ Y\in \binom{S}{r}:\left\vert Y\cap X_{1}\right\vert =r-1=\left\vert
	Y\cap X_{2}\right\vert \right\}$$
	Claim 
	$$A=\left\{ \left( X_{1}\cap X_{2}\right) \cup \{x,u\},\left(
	X_{1}\cap X_{2}\right) \cup \{y,u\},\left( X_{1}\cap X_{2}\right) \cup
	\{x,v\},\left( X_{1}\cap X_{2}\right) \cup \{y,v\}\right\}$$
	In order to prove the claim we will show that any 
	$Y \in A$ has to contain $X_{1}\cap X_{2}$. Otherwise, if there exists $z \in (X_{1}\cap X_{2}) \setminus Y$ then $Y \cap  X_{1} = X_{1}/\{z\}$ and $Y \cap X_{2} = X_{2}/\{z\}$. Therefore $Y=(X_{1}/\{z\}) \cup \{a\}$ and
	$$Y \cap X_{2}=[(X_{1}/\{z\}) \cup \{a\}]\cap X_2 =[(X_{1}/\{z\})\cap X_2] \sqcup [\{a\}\cap X_2]=[(X_{1}\cap X_2)/\{z\} ] \sqcup [\{a\} \cap X_2].$$
	Since $|Y \cap X_{2}|=r-1$, $|(X_{1}\cap X_2)/\{z\}|=r-3$ and $|\{a\}  \cap X_2|=1$ we get a contradiction. Hence, from the following three
	statements:
	\begin{itemize}
		\item 	$X_1 \cap X_2 \subset Y$
		\item  $|X_1 \cap Y|=r-1$ and  $|X_1 \cap Y|=r-1$
		\item  $X_{1}=\left( X_{1}\cap X_{2}\right) \cup \{x,y\}$ and $X_{2}=\left(X_{1}\cap X_{2}\right) \cup \{u,v\}$
	\end{itemize}  
	we get $$A=\left\{ \left( X_{1}\cap X_{2}\right) \cup \{x,u\},\left(
	X_{1}\cap X_{2}\right) \cup \{y,u\},\left( X_{1}\cap X_{2}\right) \cup
	\{x,v\},\left( X_{1}\cap X_{2}\right) \cup \{y,v\}\right\}$$
	concluding the lemma.
\end{proof}\\

\noindent
The previous lemma automatically yields
\begin{corollary}
	\label{c4}
	Given $X_1, X_2 \in \binom{S}{r}$. If $|X_1 \cap X_2|=r-2$, then 
	$$\left\vert \left\{ Y\in \binom{S}{r}:\left\vert Y\cap X_{1}\right\vert =r-1\text{ and }\left\vert Y\cap X_{2}\right\vert=r-1\right\} \right\vert =4.$$
\end{corollary}	
A straightforward conclusion of lemma \ref{intersect} is
\begin{corollary}
	\label{c5}
	Given $\{X_1, X_2,\cdots,X_m\} \subset \binom{S}{r}$ satisfying $(**)$. If
	\begin{equation}
	\label{111}
    \left\{ Y\in \binom{S}{r}:\left\vert Y\cap
	X_{i}\right\vert =r-1,\;\forall i=1,...,m\right\} \neq \emptyset
    \end{equation}  
	then
	$$\left\vert X_{i}\cap X_{j}\right\vert =r-2,\;i \neq j.$$
\end{corollary}	

The next lemma gives an upper bound on the cardinality of families of sets satisfying $(**)$ and equation \eqref{111}.
	
\begin{lemma}
	\label{17}
	If  $\{ X_{1},...,X_{m} \} \subset \binom{S}{r}$ satisfies 	
	\begin{equation}
		\label{115}
	    \left\{ Y\in \binom{S}{r}:\left\vert Y\cap
		X_{i}\right\vert =r-1,\;\forall i=1,...,m\right\} \neq \emptyset,
	    \end{equation}  
	    then $m \le r$.
\end{lemma}
\begin{proof} By hypothesis there  exists $T$ 
satisfying
$$\left\vert X_{i}\cap T\right\vert =r-1,\, i=1,...,m.$$
Then 
\begin{equation}
\label{122}
T=\left( T\cap X_{i}\right) \sqcup a_{i},\; i=1,...,m.
\end{equation}
The equation \eqref{122} together with Corollary \ref{c5} yields 
\begin{equation}
\label{123}
a_{j} \in T\cap X_{i},\; \forall j \in \{1,...,m\}\setminus \{i\}.
\end{equation}
Otherwise, there exists $a_j \in T$ and $a_j \notin T\cap X_j$ such that $a_j \notin  T\cap X_{i}$, therefore $a_j =a_i$ by \eqref{122}. This fact together with  \eqref{122} gives $T\cap X_{i}=T\cap X_{j} \subset X_i \cap X_j$ implying that $|X_i \cap X_j|=|X_i \cap T|=r-1$ which is a contradiction with the conclusion of Corollary \ref{c5}. Consequently
\begin{equation}
\label{126}
\{a_{1},...,a_{m}\}\setminus \{a_i\} \subset  T\cap X_{i}, \forall i=1,...,m.
\end{equation}
On the other hand 
$$a_{i} \neq a_{j},\; i \neq j \in \{1,...,m\},$$ 
otherwise from \eqref{123} we get $a_{i} = a_{j} \in T\cap X_{i}$, which is false by \eqref{122}. 
Therefore
$$m-1 = |\{a_{1},\cdots ,a_{m}\}\setminus \{a_i\}| \le |T\cap X_{i}|=r-1,$$
concluding that $m\le r$.
\end{proof}\\

\noindent
Let $t\leq \alpha $ and $\left\{ X_{1},...,X_{t}\right\}$ be a set satisfying the condition $(**)$. Let us define 
$$A_{i}:=\left\{ Y\in \binom{S}{r};\left\vert X_{i}\cap
Y\right\vert =r-1\right\}.$$
Then
$$\left\{ Y\in \binom{S}{r};\left\vert X_{i}\cap Y\right\vert \leq r-2,\; i=1,...,t\right\} = \binom{S}{r}\backslash \left[ \left\{
X_{1},...,X_{t}\right\} \cup \bigcup\limits_{i=1}^{t}A_{i}\right].$$

\begin{remark}
 $$\forall j=2,...,r,\; \forall\; i_{1}<\cdots <i_{j},\; \ 0\leq
\left\vert A_{i_{1}}\cap \cdots \cap A_{i_{j}}\right\vert \leq 4.$$
\end{remark}
By lemma \ref{17} the following remark holds
\begin{remark}
 $$\forall
j=r+1,...,\alpha,\; \forall\; i_{1}<\cdots <i_{j}, \; A_{i_{1}}\cap \cdots \cap A_{i_{j}}=\emptyset. $$
\end{remark} 
From the previous remark we have the following lemmas

\begin{lemma}
\label{l20}
Let	$\left\{ X_{1},...,X_{t}\right\}$ be a set satisfying the condition $(**)$. If $t\leq r$,  then
$$tr(n-r)-4\sum\limits_{\overset{j\text{even}}{j\in \{2,...,t\}}}\binom{t
}{j}\le \left\vert \bigcup\limits_{i=1}^{t}A_{i}\right\vert \leq
tr(n-r)+4\sum\limits_{\overset{j\text{ odd}}{j\in \{2,...,t\}}}\binom{t}{j}.$$
\end{lemma}
\begin{proof}
Right hand	
\begin{eqnarray} 
\left\vert \bigcup\limits_{i=1}^{t}A_{i}\right\vert
&=&\sum\limits_{i=1}^{t}\left\vert A_{i}\right\vert
+\sum\limits_{j=2}^{t}\sum\limits_{\underset{i_{1},...,i_{j}\in \{1,...,t\}}{
i_{1}<\cdots <i_{j}}}\left( -1\right) ^{j-1}\left\vert A_{i_{1}}\cap \cdots \cap A_{i_{j}}\right\vert \\
&\leq& tr(n-r)+\sum\limits_{\overset{j\text{odd}}{j\in \{2,...,t\}}%
}\sum\limits_{\underset{i_{1},...,i_{j}\in \{1,...,t\}}{i_{1}<\cdots <i_{j}}}\left\vert A_{i_{1}}\cap \cdots \cap A_{i_{j}}\right\vert \\
&\leq&
tr(n-r)+4\sum\limits_{\overset{j\text{odd}}{j\in \{2,...,t\}}}\binom{t}{j}. 
\end{eqnarray}
Left hand
\noindent
\begin{eqnarray}
\left\vert \bigcup\limits_{i=1}^{t}A_{i}\right\vert
&=&\sum\limits_{i=1}^{t}\left\vert A_{i}\right\vert
+\sum\limits_{j=2}^{t}\sum\limits_{\underset{i_{1},...,i_{j}\in \{1,...,t\}}{i_{1}<\cdots <i_{j}}}\left( -1\right) ^{j-1}\left\vert A_{i_{1}}\cap \cdots \cap A_{i_{j}}\right\vert\\
&\geq& tr(n-r)-\sum\limits_{\overset{j\text{even}}{j\in \{2,...,t\}}%
}\sum\limits_{\underset{i_{1},...,i_{j}\in \{1,...,t\}}{i_{1}<\cdots <i_{j}}}\left\vert A_{i_{1}}\cap \cdots \cap A_{i_{j}}\right\vert \\
&\geq& tr(n-r)-4\sum\limits_{\overset{j\text{even}}{j\in \{2,...,t\}}}\binom{t
}{j}.
\end{eqnarray}
\end{proof}\\
\noindent
A similar lemma for $t\geq r$ is presented
\begin{lemma}
\label{l21}
If $t\geq r$, then
$$
tr(n-r)-4\sum\limits_{\overset{j\text{even}}{j\in \{2,...,r\}}}\binom{t%
}{j}\le \left\vert \bigcup\limits_{i=1}^{t}A_{i}\right\vert \le tr(n-r)+4\sum\limits_{\overset{j\text{odd}}{j\in \{2,...,r\}}}\binom{t}{j}.
$$
\end{lemma}
\begin{proof}
Right hand
\begin{eqnarray}
\left\vert \bigcup\limits_{i=1}^{t}A_{i}\right\vert &=& \sum\limits_{i=1}^{t}\left\vert A_{i}\right\vert
+\sum\limits_{j=2}^{t}\sum\limits_{\underset{i_{1},...,i_{j}\in \{1,...,t\}}{i_{1}<\cdots <i_{j}}}\left( -1\right)^{j-1}\left\vert A_{i_{1}}\cap \cdots \cap A_{i_{j}}\right\vert\\
&=&\sum\limits_{i=1}^{t}\left\vert A_{i}\right\vert
+\sum\limits_{j=2}^{r}\sum\limits_{\underset{i_{1},...,i_{j}\in \{1,...,t\}}{i_{1}<\cdots <i_{j}}}\left( -1\right) ^{j-1}\left\vert A_{i_{1}}\cap \cdots \cap A_{i_{j}}\right\vert \\
&\leq& tr(n-r)+\sum\limits_{\overset{j\text{odd}}{j\in \{2,...,r\}}
}\sum\limits_{\underset{i_{1},...,i_{j}\in \{1,...,t\}}{i_{1}<\cdots <i_{j}}}\left\vert A_{i_{1}}\cap \cdots \cap A_{i_{j}}\right\vert \\
&\leq & tr(n-r)+4\sum\limits_{\overset{j\text{odd}}{j\in \{2,...,r\}}}\binom{t}{j}.
\end{eqnarray}
Left hand
\noindent
\begin{eqnarray}
\left\vert \bigcup\limits_{i=1}^{t}A_{i}\right\vert
&=&\sum\limits_{i=1}^{t}\left\vert A_{i}\right\vert
+\sum\limits_{j=2}^{r}\sum\limits_{\underset{i_{1},...,i_{j}\in \{1,...,r\}}{i_{1}<\cdots <i_{j}}}\left( -1\right) ^{j-1}\left\vert A_{i_{1}}\cap \cdots \cap A_{i_{j}}\right\vert\\
&\geq& tr(n-r)-\sum\limits_{\overset{j\text{even}}{j\in \{2,...,r\}}%
}\sum\limits_{\underset{i_{1},...,i_{j}\in \{1,...,r\}}{i_{1}<\cdots <i_{j}}}\left\vert A_{i_{1}}\cap \cdots \cap A_{i_{j}}\right\vert \\
&\geq& tr(n-r)-4\sum\limits_{\overset{j\text{even}}{j\in \{2,...,r\}}}\binom{t
}{j}.
\end{eqnarray}
\end{proof}\\
\noindent
An immediate consequence of lemmas \ref{l20} and \ref{l21} is the following corollary. 
\begin{corollary} 
\label{cor22}
If $t\le r$, then 
\begin{eqnarray}
\binom{n}{r}-tr(n-r)-4\sum\limits_{\overset{j\text{even}}{j\in \{2,...,t\}}}\binom{t}{j} & \leq & \left\vert \left\{ Y\in \binom{S}{r};\left\vert X_{i}\cap Y\right\vert \leq r-2, \forall i=1,...,t\right\} \right \vert \nonumber \\
&=&\left\vert \binom{S}{r}\backslash \bigcup\limits_{i=1}^{t}A_{i}\right\vert=\binom{n}{r}-\left\vert \bigcup\limits_{i=1}^{t}A_{i}\right\vert \nonumber \\ 
&\leq& \binom{n}{r}-tr(n-r)+4\sum\limits_{\overset{j\text{odd}}{j\in \{2,...,t\}}}\binom{t}{j}.
\end{eqnarray}
If $t>r$, then 
\begin{eqnarray}
\binom{n}{r}-tr(n-r)-4\sum\limits_{\overset{j\text{even}}{j\in \{2,...,r\}}}\binom{t}{j} & \leq & \left\vert \left\{ Y\in \binom{S}{r};\left\vert X_{i}\cap Y\right\vert \leq r-2, \forall i=1,...,t\right\} \right \vert \nonumber \\
&=&\left\vert \binom{S}{r}\backslash \bigcup\limits_{i=1}^{t}A_{i}\right\vert=\binom{n}{r}-\left\vert \bigcup\limits_{i=1}^{t}A_{i}\right\vert \nonumber \\ 
&\leq& \binom{n}{r}-tr(n-r)+4\sum\limits_{\overset{j\text{odd}}{j\in \{2,...,r\}}}\binom{t}{j}.
\end{eqnarray}
\end{corollary}
\noindent 
The previous corollary leads to a theorem that counts the number of subsets of $\binom{S}{r}$ satisfying the condition $(**)$. 
\begin{theorem}
	
If $t\leq r$, then 

$$\left\vert \left\{ \left\{ X_{1},...,X_{t}\right\} ;\text{satisfying }(\ast\ast )\right\} \right\vert \leq $$

$$\leq \frac{1}{t!}\binom{n}{r}\left[ \binom{n}{r}-\left[ r(n-r)+1\right] 
\right] \prod\limits_{h=2}^{t}\left[ \binom{n}{r}-h\left[ r(n-r)+1\right]
+4\sum\limits_{\overset{j\text{odd}}{j\in \{2,...,h\}}}\binom{h}{j}\right]$$
\noindent
and

$$\left\vert \left\{ \left\{ X_{1},...,X_{t}\right\} ;\text{satisfying }(\ast
\ast )\right\} \right\vert \geq $$
$$\geq \frac{1}{t!}\binom{n}{r}\left[ \binom{n}{r}-\left[
r(n-r)+1\right] \right] \prod\limits_{h=2}^{t}\left[ \binom{n}{r}-h\left[
r(n-r)+1\right] -4\sum\limits_{\overset{j\text{even}}{j\in \{2,...,h\}}}
\binom{h}{j}\right].$$
\noindent
On the other hand, if $r<t\leq \alpha$, then

$$\left\vert \left\{ \left\{ X_{1},...,X_{t}\right\} ;\text{satisfying }(\ast\ast )\right\} \right\vert \leq $$

$$\leq \frac{1}{t!}\binom{n}{r}\left[ \binom{n}{r}-\left[
r(n-r)+1\right]  \right] \prod\limits_{h=2}^{r}\left[ \binom{n}{r}-h \left[ r(n-r)+1\right] +4\sum\limits_{\overset{j\text{odd}}{j\in \{2,...,h\}}}\binom{h}{j}\right]  \cdot $$
$$\cdot \prod\limits_{h=r+1}^{t}\left[ \binom{n}{r}-h\left[ r(n-r)+1\right] +4\sum\limits_{\overset{j\text{odd}}{j\in \{2,...,r\}}}\binom{h}{j}\right]$$
\noindent
and

$$\left\vert \left\{ \left\{ X_{1},...,X_{t}\right\} ;\text{satisfying 
}(\ast \ast )\right\} \right\vert \geq $$
$$\geq \frac{1}{t!}\binom{n}{r}\left[ \binom{n}{r}-\left[
r(n-r)+1\right] \right]  \prod\limits_{h=2}^{r}\left[ \binom{n}{r}-h \left[ r(n-r)+1\right] -4\sum\limits_{\overset{j\text{ odd}}{j\in \{2,...,h\}}}\binom{h}{j}\right] \cdot $$
$$\cdot \prod\limits_{h=r+1}^{t}\left[ \binom{n}{r}-h\left[ r(n-r)+1\right] -4\sum\limits_{\overset{j \text{ odd}}{j\in \{2,...,r\}}}\binom{h}{j}\right].$$
\end{theorem}

From this theorem we can find an upper and lower bound for the 
number of sparse-paving matroids as follows

\begin{theorem}
	The following upper and lower bounds are valid
$$\left\vert Sparse_{n,r} \right\vert \leq 1+\binom{n}{r}+ \frac{1}{2}\binom{n}{r}\left[ \binom{n}{r}%
-r(n-r)-1\right]  +$$
$$+\sum\limits_{t=2}^{r} \frac{1}{t!}\binom{n}{r}\left[ \binom{n}{r}-%
\left[ r(n-r)+1\right] \right] \prod\limits_{h=2}^{t}\left[ \binom{n}{r}-h
\left[ r(n-r)+1\right] +4\sum\limits_{\overset{j\text{even}}{j\in \{2,...,h\}%
}}\binom{h}{j}\right] +$$ 

$$+
\sum\limits_{t=r+1}^{\alpha }\frac{1}{t!}\binom{n}{r}\left[ \binom{n}{r}-%
\left[ r(n-r)+1\right] \right]  \prod\limits_{h=3}^{r}\left[ \binom{n}{%
	r}-h\left[ r(n-r)+1\right] +4\sum\limits_{\overset{j\text{even}}{j\in
		\{2,...,h\}}}\binom{h}{j}\right]  \cdot $$ 
$$
\cdot \prod\limits_{h=r+1}^{t}\left[ \binom{n}{r}-h\left[ r(n-r)+1\right] +4\sum\limits_{\overset{j\text{odd}}{j\in \{2,...,r\}}}\binom{h}{j}\right]
$$

$$\left\vert Sparse_{n,r} \right\vert \geq 1+\binom{n}{r}+ \frac{1}{2}\binom{n}{r} \binom{n}{r}-r(n-r)-1+$$
$$+\sum\limits_{t=3}^{r} \frac{1}{t!}\binom{n}{r}\left[ \binom{n}{r}-
\left[ r(n-r)+1\right] \right] \prod\limits_{h=2}^{t} \left [\binom{n}{r}-h
\left[ r(n-r)+1\right] -4\sum\limits_{\overset{j\text{odd}}{j\in
		\{2,...,h\}}}\binom{h}{j}\right] +$$

$$+
\sum\limits_{t=r+1}^{\alpha }\frac{1}{t!}\binom{n}{r}\left[ \binom{n}{r}-%
\left[ r(n-r)+1\right] \right]  \prod\limits_{h=2}^{r}\left[ \binom{n}{%
	r}-h\left[ r(n-r)+1\right] -4\sum\limits_{\overset{j\text{even}}{j\in
		\{2,...,h\}}}\binom{h}{j}\right]  \cdot $$
	$$ 
\cdot \prod\limits_{h=r+1}^{t}\left[ \binom{n}{r}-h\left[ r(n-r)+1\right] -4\sum\limits_{\overset{j\text{odd}}{j\in \{2,...,r\}}}\binom{h}{j}\right].
$$
\end{theorem}
The next section provides a different idea to construct families of sets satisfying $(**)$. 

\section{ A method to construct sets\textbf{\ }$\mathcal{U}$
with property $(\ast \ast )\;:\; \forall X,Y\in \mathcal{U}$ with $
X\neq Y,\; \left\vert X\cap Y\right\vert \leq r-2$.}
\label{4s}

In this section, we will construct matrices of $r$-subsets of $S$ from a
fixed $r$-subset $X$ having the following properties:

\begin{description}
\item{a)}  Any $r$-subset of $S$ is an entry of
exactly one of these matrices.

\item{b)} In each matrix, any two entries which are in different
rows and different columns have intersection less or equal to $r-2$.

\item{c)} Any two entries in different matrices  $S_{h}$
and $S_{h}$ with $\left\vert h-h^{\prime }\right\vert \geq 2$, 
have intersection with cardinality less than or equal to $r-2$.
\end{description}

Let $S$ be a set of cardinality $n$, $2\leq r\leq n-1$. Fix $X\in \binom{S}{r}$. For each $0\leq h\leq r$, let 
$$\binom{X}{h}
=\left \{A_{1}^{(h)},...,A_{\binom{r}{h}}^{(h)} \right \}$$ 
be the $h$-subsets of $X$ and $$\binom{S\backslash X}{r-h}=\left \{Z_{1}^{(h)},...,Z_{\binom{n-r}{r-h}}^{(h)} \right \}$$
be the $(r-h)$-subsets of $S\backslash X$. For each $0\leq h\leq r$ and $
n-r\geq r-h$, we build the following $\binom{\left\vert S\backslash
X\right\vert }{r-h}\times \binom{\left\vert X\right\vert }{h}$-matrix $
{\LARGE s}_{h}$ as follow

$$
{\LARGE s}_{h}:=\left[ 
\begin{array}{ccc}
A_{1}^{(h)}\cup Z_{1}^{(h)} & \cdots & A_{\binom{r}{h}}^{(h)}\cup Z_{1}^{(h)}
\\ 
A_{1}^{(h)}\cup Z_{2}^{(h)} & \cdots & A_{\binom{r}{h}}^{(h)}\cup Z_{2}^{(h)}
\\ 
\vdots & \vdots & \vdots \\ 
A_{1}^{(h)}\cup Z_{\binom{n-r}{r-h}}^{(h)} & \cdots & A_{\binom{r}{h}
}^{(h)}\cup Z_{\binom{n-r}{r-h}}^{(h)}
\end{array}
\right] _{\binom{S\backslash X}{r-h}\times \binom{X}{h}}.
$$
From the previous construction we deduce the following lemma

\begin{lemma}
\label{Construction} The matrices ${\LARGE s}_{h}$ have the following
properties
\begin{description}
\item{a)} $\forall\, Y \in \binom{S}{r},$ there exists a unique $
0\leq h\leq r$ such that $Y$ is an entry of $S_{h}$.

\item{b)} Let  $1\leq i\neq j\leq \binom{n-r}{r-h}$ and $1\leq
t\neq k\leq \binom{r}{h}$. Then $$\left\vert \left( A_{t}^{(h)}\cup
Z_{i}^{(h)}\right) \cap \left( A_{k}^{(h)}\cup Z_{j}^{(h)}\right)
\right\vert \leq r-2.$$
That is, any pair of entries in different
columns and different rows have intersection smaller than $r-2$.

\item{c)} For $0\leq h,h^{\prime }\leq r$ such that
 $\left\vert h-h^{\prime }\right\vert \geq 2$ and for all $i,j,t$ 
and $k$, 
$$\left\vert \left( A_{t}^{(h)}\cup Z_{i}^{(h)}\right) \cap
\left( A_{k}^{(h^{\prime })}\cup Z_{j}^{(h^{\prime })}\right) \right\vert
\leq r-2.$$
\end{description}
\end{lemma}

\begin{proof} The item a) follows directly from the definition
of the matrices $S_{h}$. In order to prove the item b) we observe that if 
$A_{t}^{(h)}\neq A_{k}^{(h)}\subseteq X$ and $Z_{i}^{(h)}\neq
Z_{j}^{(h)}\subseteq S\backslash X$, then 

\begin{eqnarray*}
\left\vert \left( A_{t}^{(h)}\cup Z_{i}^{(h)}\right) \cap \left(
A_{k}^{(h)}\cup Z_{j}^{(h)}\right) \right\vert &=&\left\vert A_{t}^{(h)}\cap
A_{k}^{(h)}\right\vert +\left\vert Z_{i}^{(h)}\cap Z_{j}^{(h)}\right\vert \\
&\leq & \left( h-1\right) +\left( r-h-1\right) =r-2.
\end{eqnarray*}

In regard to item c) we need to prove that 
$$\left\vert \left(
A_{t}^{(h)}\cup Z_{i}^{(h)}\right) \cap \left( A_{k}^{(h^{\prime })}\cup
Z_{j}^{(h^{\prime })}\right) \right\vert \leq r-2$$ 
for $0\leq h,h^{\prime}\leq r$ with $\left\vert h-h^{\prime }\right\vert \geq 2$ and $1\leq
i,j\leq \binom{n-r}{r-h}$, $1\leq t,k\leq \binom{r}{h}$.  Since $\left\vert
h-h^{\prime }\right\vert \geq 2$, we can assume that $h=h^{\prime }+m$ with 
$2\leq m\leq r-h^{\prime }$. Then 
\begin{eqnarray*}
\left\vert \left( A_{t}^{(h)}\cup Z_{i}^{(h)}\right) \cap \left(
A_{k}^{(h^{\prime })}\cup Z_{j}^{(h^{\prime })}\right) \right\vert
&=&\left\vert A_{t}^{(h)}\cap A_{k}^{(h^{\prime })}\right\vert +\left\vert
Z_{i}^{(h)}\cap Z_{j}^{(h^{\prime })}\right\vert \\
&\leq & h+(r-h^{\prime
})=r-m\leq r-2.
\end{eqnarray*}
\end{proof}

\noindent
Now we will give some examples of the matrices ${\LARGE s}_{h}$. \\

\noindent
\textbf{Example 1:} Let $S=\{1,2,3,4,5,6\},$ $r=3$ and fix $X=\{1,2,3\}$.

$$
{\LARGE s}_{0}:=\left[ 
\begin{array}{c}
\{4,5,6\}
\end{array}
\right] _{\binom{S\backslash X}{3}\times \binom{X}{0}},$$

$${\LARGE s}_{1}:=\left[ 
\begin{array}{ccc}
\{1\}\cup \{4,5\} & \{2\}\cup \{4,5\} & \{3\}\cup \{4,5\} \\ 
\{1\}\cup \{4,6\} & \{2\}\cup \{4,6\} & \{3\}\cup \{4,6\} \\ 
\{1\}\cup \{5,6\} & \{2\}\cup \{5,6\} & \{3\}\cup \{5,6\}
\end{array}
\right] _{\binom{S\backslash X}{2}\times \binom{X}{1}},$$

$${\LARGE s}_{2}:=\left[ 
\begin{array}{ccc}
\{1,2\}\cup \{4\} & \{1,3\}\cup \{4\} & \{2,3\}\cup \{4\} \\ 
\{1,2\}\cup \{5\} & \{1,3\}\cup \{5\} & \{2,3\}\cup \{5\} \\ 
\{1,2\}\cup \{6\} & \{1,3\}\cup \{6\} & \{2,3\}\cup \{6\}%
\end{array}
\right] _{\binom{S\backslash X}{1}\times \binom{X}{2}}$$
and

$${\LARGE s}_{3}:=\left[ \{1,2,3\}\right] _{\binom{S\backslash X}{0}\times 
\binom{X}{3}}.$$
\noindent
 In the next subsection using the matrices ${\LARGE s}_{h}$ we will
find partitions of the set $\binom{S}{r}$ satisfying property $(\ast \ast )$.

\subsection{A partition of $\binom{S}{r}$ by
subsets satisfying property $(\ast \ast )$ }
Now we will construct  $\mathcal{U}^{\prime }$s having property $(\ast \ast
)$ which form a partition of $\binom{S}{r}$. Let $S=\{1,2,...,n\}$, $2\leq
r\leq n-1$ and let $X\in \binom{S}{r}$ be fixed. Let $0\leq h\leq r$ and
take the$\binom{\left\vert S\backslash X\right\vert }{r-h}\times \binom{
\left\vert X\right\vert }{h}$-matrix ${\LARGE s}_{h}$. By the previous lemma
item b), we can make $\max \{\binom{n-r}{r-h},\binom{r}{h}\}$ different sets
consisting of the entries of $S_{h}$ satisfying property $(\ast \ast )$.
Namely, take each set  with the entries of each major diagonal of $S_{h}$.
In this way, we get $\max \{\binom{n-r}{r-h},\binom{r}{h}\}$ different sets
of cardinality $\min \{\binom{n-r}{r-h},\binom{r}{h}\}$. Graphically
speaking, from
$$\left [ 
\begin{array}{ccc}
\bullet _{1} & \triangleright _{2} & \circ _{3} \\ 
\ast _{1} & \bullet _{2} & \triangleright _{3} \\ 
\circ _{1} & \ast _{2} & \bullet _{3} \\ 
\triangleright _{1} & \circ _{2} & \ast _{3}
\end{array}
\right ] _{4\times 3} \longleftrightarrow \left [
\begin{array}{ccc}
\bullet _{1} &  &  \\ 
\ast _{1} & \bullet _{2} &  \\ 
\circ _{1} & \ast _{2} & \bullet _{3} \\ 
\triangleright _{1} & \circ _{2} & \ast _{3} \\ 
& \triangleright _{2} & \circ _{3} \\ 
&  & \triangleright _{3}
\end{array}
\right ], $$

\noindent 
we obtain $\{\bullet _{1},\bullet _{2},\bullet _{3}\},\{\ast
_{1},\ast _{2},\ast _{3}\},\{\circ _{1},\circ _{2},\circ _{3}\}$ and $
\{\triangleright _{1},\triangleright _{2},\triangleright _{3}\}$. In the
following we will show that the set constructed in this way fulfill the
property $(\ast \ast )$. We start by giving some definitions

$$
\begin{array}{ccc}
\sigma _{j}^{(h)}: & \left\{ 1,2,...,\min \left\{ \binom{n-r}{r-h},\binom{r}{h}
\right\} \right\}   & \rightarrow   \left\{ 1,2,...,\max \left\{ \binom{n-r}{%
r-h},\binom{r}{h}\right\} \right\}  \\ 
& t  & \mapsto   \left[ j+t-1\right] _{\text{
mod}\max \left\{ \binom{n-r}{r-h},\binom{r}{h}\right\} }%
\end{array}.
$$
For all\textbf{\ }$0\leq h\leq r$ and $j=1,...,\max \left\{ \binom{n-r}{%
r-h},\binom{r}{h}\right\} $ we define the sets

$${\LARGE s}_{h}(j):=\bigsqcup\limits_{t=1}^{\min \left\{ \binom{n-r}{r-k},%
\binom{r}{k}\right\} }\left\{ A_{t}^{(h)}\cup Z_{\sigma
_{j}^{(h)}(t)}^{(h)}\right\} .$$

\begin{lemma} The sets ${\LARGE s}_{h}(j):=\bigsqcup\limits_{t=1}^{\min
\left\{ \binom{n-r}{r-k},\binom{r}{k}\right\} }\left\{ A_{t}^{(h)}\cup
Z_{\sigma _{j}^{(h)}(t)}^{(h)}\right\} $ satisfy the property $(\ast \ast )$.
\end{lemma}

\begin{proof} Without loss of generality we consider the case $%
\binom{n-h}{r-h}\geq \binom{r}{h}$. Then $\max \{\binom{n-r}{r-h},\binom{r}{h%
}\}=\binom{n-h}{r-h}.$ In this case, the function $\sigma _{j}^{(h)}$ take
the form

$$\sigma _{j}^{(h)}:\{1,2,...,\binom{r}{h}\}\rightarrow \{1,2,...,\binom{n-r}{
r-h}\},\; \sigma _{j}^{(h)}(t)=\left[ j+t-1\right] _{\text{mod}\binom{n-r}{
r-h}}$$
which is an injective function and each pair $(t,$ $\sigma _{j}^{(h)})$
provides a different element of the $j$-th diagonal of ${\LARGE s}_{h}$. By
item b) of the  lemma \ref{Construction} two different elements  $A_{t}^{(h)}\cup Z_{\sigma
_{j}^{(h)}(t)}^{(h)}$ and $A_{t}^{(h)}\cup Z_{\sigma _{j}^{(h)}(t)}^{(h)}$ have intersection in at most $r-2$ elements, consequently

$${\LARGE s}_{h}(j):=\bigsqcup\limits_{t=1}^{\binom{r}{h}}\left\{
A_{t}^{(h)}\cup Z_{\sigma _{j}^{(h)}(t)}^{(h)}\right\} $$ 
satisfies the property $(\ast \ast )$.
\end{proof}\\

Now,  we will give bigger subsets of $\binom{S}{r}$ which still
have property $(\ast \ast )$.

\begin{lemma} 
The sets
$$\mathcal{U}_{j}^{(odd)}:=\bigsqcup\limits_{0\leq h\leq r\text{, }h\text{
odd}}{\LARGE s}_{h}(j)$$ 
for $$j=1,...,\underset{0\leq k\leq
r,\text{ }k\text{ odd}}{\max }\left\{ \max \left\{ \binom{n-r}{r-k},\binom{r%
}{k}\right\} \right\} $$
and 
$$\mathcal{U}_{j}^{(even)}:=\bigsqcup\limits_{0\leq h\leq r\text{, }h\text{
even}}{\LARGE s}_{h}(j)$$ 
for 
$$j=1,...,\underset{0\leq k\leq r,
\text{ }k\text{ even}}{\max }\left\{ \max \left\{ \binom{n-r}{r-k},\binom{r}{
k}\right\} \right\} ,$$
where in both cases, for all $h$, $S_{h}(j):=\emptyset $ if $
j>\max \left\{ \binom{n-r}{r-h},\binom{r}{h}\right\} $ satisfies condition $
(\ast \ast )$.
\end{lemma}

\begin{proof}
Let fix $j$. Take $A,B\in \mathcal{U}_{j}^{(odd)}$ then
we have two cases:

\begin{enumerate}
\item $A,B\in {\LARGE s}_{h}(j)$ for some $h$.

\item $A\in {\LARGE s}_{h}(i)$ and $B\in {\LARGE s}_{h'}(j)$ for  $h \neq  h'$.
\end{enumerate}

The case 1 follows from lemma \ref{Construction} item b). In 
case 2  we have that  $\left\vert h-h^{\prime }\right\vert \geq 2$ and the
conclusion follows from lemma \ref{Construction} c). 
\end{proof}\\

\noindent
One readily see that 

$$\left\{ \mathcal{U}_{j}^{(odd)}\right\} _{j=1,...,\underset{0\leq k\leq r,%
\text{ }k\text{ odd}}{\max }\binom{n-r}{r-k}}\bigcup \left\{ \mathcal{U}%
_{j}^{(even)}\right\} _{j=1,...,\underset{0\leq k\leq r,\text{ }k\text{ even}%
}{\max }\binom{n-r}{r-k}}$$
is a partition of $\binom{S}{r}$ with
$$\gamma :=\max_{0\leq h\leq r\text{, }h\text{ odd}} \left \{\max \left \{\binom{n-r}{r-h},
\binom{r}{h} \right \}+\max_{0\leq h\leq r\text{, }h\text{ even}} \left \{\max \left \{\binom{n-r}{r-h},\binom{r}{h} \right \} \right \}\right \}$$
elements. \\

Now we will give an example in order to facilitate the understanding of the above  construction. \\

\textbf{Example 2. }From example 1, we have: $S=\{1,...,6\}$ and $r=3.$
$${\LARGE s}_{1}:=\left[ 
\begin{array}{ccc}
\text{\textbf{(a)}}\{1,4,5\} & \text{\textbf{(c)}}\{2,4,5\} & \text{\textbf{
(b)}}\{3,4,5\} \\ 
\text{\textbf{(b)}}\{1,4,6\} & \text{\textbf{(a)}}\{2,4,6\} & \text{\textbf{
(c)}}\{3,4,6\} \\ 
\text{\textbf{(c)}}\{1,5,6\} & \text{\textbf{(b)}}\{2,5,6\} & \text{\textbf{
(a)}}\{3,5,6\}
\end{array}
\right] .$$
By lemma \ref{Construction} b), each of $\left\{ \{1,4,5\},\{2,4,6\},\{3,5,6\}\right\} ,
\, \left\{ \{1,4,6\},\{2,5,6\},\{3,4,5\}\right\} $ and $\left\{
\{1,5,6\},\{2,4,5\},\{3,4,6\}\right\} $ satisfies property $(\ast \ast ).$
And by lemma \ref{Construction} c),

\begin{eqnarray*}
\mathcal{U}_{1}^{(even)}& =&\left\{
\{4,5,6\},\{1,2,4\},\{1,3,5\},\{2,3,6\}\right\} , \\ 
\mathcal{U}_{2}^{(even)}&=&\left\{ \{1,2,5\},\{1,3,6\},\{2,3,4\}\right\} , \\ 
\mathcal{U}_{3}^{(even)}&=&\left\{ \{1,2,6\},\{1,3,4\},\{2,3,5\}\right\} , \\ 
\mathcal{U}_{1}^{(odd)}&=&\left\{ \{1,4,5\},\{2,4,6\},\{3,5,6\}\right\} , \\ 
\mathcal{U}_{2}^{(odd)}&=&\left\{ \{1,4,6\},\{2,5,6\},\{3,4,5\}\right\} \\
\text{and} && \\ 
\mathcal{U}_{3}^{(odd)}&=&\left\{ \{1,5,6\},\{2,4,5\},\{3,4,6\}\right\} .
\end{eqnarray*}
  
That is, each $\mathcal{U}_{i}$ satisfies property $(\ast \ast )$ and $%
\binom{S}{3}=\sqcup _{i=1}^{3}\mathcal{U}_{i}^{(odd)}\sqcup \sqcup _{i=1}^{3}
\mathcal{U}_{i}^{(even)}$.

%\subsection{A lower bound on the number of $\left\vert
%Sparse_{n,r}\right\vert$ }
%\begin{lemma}
%\label{lowerbound}
%Given a set $S$ of $n$ elements the number $\left\vert
%Sparse_{n,r}\right\vert$ can be bounded from below as
%$$
%\left\vert Sparse_{n,r}\right\vert \ge \max \left \{ \prod_{h\; \text{even} } \left ( \sum_{\mathcal{C}_{h}^{r}  \subset \mathcal{P}(X) }\sum_{A \in \mathcal{C}_{h}^{r} }  \binom{n-r}{r-h}^{|A|} \right ), \prod_{h\; \text{odd} } \left ( \sum_{\mathcal{C}_{h}^{r}  \subset \mathcal{P}(X) }\sum_{B \in \mathcal{C}_{h}^{r} }   \binom{n-r}{r-h}^{|B|} \right ) \right \}
%$$
%where  $\mathcal{C}_{h}^{r} \subset \mathcal{P}(X)$ and each $A \in \mathcal{C}_{h}^{r}$ satisfies $(**)$. Also 
%
%$$
%\left\vert Sparse_{n,r}\right\vert \ge \max \left \{ \prod_{h\; \text{even}} \left ( \sum_{\mathcal{C}_{h}^{n-r}  \subset \mathcal{P}(S \setminus X) }\sum_{A \in \mathcal{C}_{h}^{n-r} }  \binom{r}{h}^{|A|} \right ), \prod_{h\; \text{odd}} \left ( \sum_{\mathcal{C}_{h}^{n-r}  \subset \mathcal{P}(S\setminus X) }\sum_{B \in \mathcal{C}_{h}^{n-r} }   \binom{r}{h}^{|B|} \right ) \right \}
%$$
%where  $\mathcal{C}_{h}^{n-r} \subset \mathcal{P}(X)$ and each $A \in \mathcal{C}_{h}^{n-r}$ satisfies $(**)$. 
%\end{lemma}
%\begin{proof}
%To be written
%\end{proof}

%These lower bounds can be stated in a more concrete form as the next lemma shows

$(\ast )$

\textbf{takane@matem.unam.mx (email contact)}

Instituto de Matem\'{a}ticas

Universidad Nacional Aut\'{o}noma de M\'{e}xico (UNAM)

www.matem.unam.mx

$(\ast \ast )$

\textbf{boris.mederos@uacj.mx}

\textbf{gtapia@uacj.mx}

Instituto de Ingenier\'{\i}a y Tecnolog\'{\i}a

Universidad Aut\'{o}noma de Ciudad Ju\'{a}rez

http://www.uacj.mx/IIT

$(\ast \ast \ast )$

\textbf{bzs@unam.mx}

Facultad de Ciencias

Universidad Aut\'{o}noma del Estado de M\'{e}xico

http://www.uaemex.mx/fciencias/

\end{document}